\def\real{{\tt I\kern-.2em{R}}}
\def\nat{{\tt I\kern-.2em{N}}}

\def\realp#1{{\tt I\kern-.2em{R}}^#1}
\def\natp#1{{\tt I\kern-.2em{N}}^#1}
\def\hyper#1{\ ^*\kern-.2em{#1}}

\def\hyperrealp#1{{\tt ^*{I\kern-.2em{R}}}^#1} 
\def\hypernat{{^*{\nat }}}
\def\hypernatp#1{{{^*{{\tt I\kern-.2em{N}}}}}^#1} 
\def\eskip{\hskip.25em\relax}

\def\Hyper#1{\hyper {\eskip #1}}
\def\leaderfill{\leaders\hbox to 1em{\hss.\hss}\hfill}
\def\srealp#1{{\rm I\kern-.2em{R}}^#1}

\def\power#1{{{\cal P}(#1)}}

\def\pars{\par\smallskip}
\def\parm{\par\medskip}
\def\r#1{{\rm #1}}
\def\b#1{{\bf #1}}
\def\ref#1{$^{#1}$}

\def\sig{{^\sigma}}
\def\m@th{\mathsurround=0pt}
\def\rightarrowfill{$\m@th \mathord- \mkern-6mu \cleaders\hbox{$\mkern-2mu 
\mathord- \mkern-2mu$}\hfil \mkern-6mu \mathord\rightarrow$}
\def\leftarrowfill{$\mathord\leftarrow
\mkern -6mu \m@th \mathord- \mkern-6mu \cleaders\hbox{$\mkern-2mu 
\mathord- \mkern-2mu$}\hfil $}
\def\noarrowfill{$\m@th \mathord- \mkern-6mu \cleaders\hbox{$\mkern-2mu 
\mathord- \mkern-2mu$}\hfil$}
\def\orgate{$\bigcirc \kern-.80em \lor$}
\def\andgate{$\bigcirc \kern-.80em \land$}
\def\inverter{$\bigcirc \kern-.80em \neg$}

\magnification=\magstep0 
\tolerance 10000
\baselineskip  14pt
\hoffset=.25in
\hsize 6.00 true in

\font\ninerm=cmr9

\vsize 8.75 true in
{\quad}\parm
The following is a direct copy of this published paper$,$ but with numerously many printer's errors corrected and with a few additional remarks.\par\medskip
\hrule\smallskip
\hrule\smallskip
\noindent {\bf Herrmann$,$ R. A.\par
\noindent Kobe J. Math.}\par
\noindent \b 4(1987) 1-14\par
\bigskip
\centerline{\bf NONSTANDARD CONSEQUENCE OPERATORS}\bigskip
\centerline{By Robert A. Herrmann}\medskip
\centerline{\ninerm (Dedicated to Professor K. Is\'eki)}
\centerline{\ninerm (Receiver October 22, 1984)}\bigskip
\indent {\bf 1. Introduction}\parm
In 1963, Abraham Robinson applied his newly discovered nonstandard analysis to formal first-order languages and developed a nonstandard logic [11] relative to the ``truth'' concept and structures. Since that time not a great deal of fundamental research has been attempted in this specific area with one notable exception [3]. However, when results from this discipline are utilized they have yielded some highly significant and important developments such as those obtained by Henson [4]. \pars
The major purpose for this present investigation is to institute formally a more general study than previously pursued. In particular, we study nonstandard logics relative to consequence operators [2] [6] [12] [13] defined on a nonstandard language. Since the languages considered are not obtained by the usual constructive methods, then this will necessitate the construction of an entirely new foundation distinctly different from Robinson's basic embedding techniques. Some very basic results of this research were very briefly announced in a previous report [6].\pars
In order to remove ambiguity from the definition of the ``finite'' consequence operator, the definition of ``finite'' is the ordinary definition in that the empty set is finite and any nonempty set $\r A$ is finite if and only if there exists a bijection $\rm f\colon A \to [1,n],$ where $\rm [1,n] = \{x\mid n\in \nat,\ 1\leq x \leq n\}$  ($\nat$ is the set of natural numbers with zero). Unless otherwise stated, all sets $\r B$ that are infinite will also be assumed to be Dedekind-infinite. This occurs when a set $\r B$ is denumerable, since $\r B$ inherits a well-ordering from $\nat$, or $\r B$ is well-ordered [2, p. 248], or the Axiom of Choice is assumed. We note that within mathematics one is always allowed to make a finite choice from finitely many nonempty sets, among others [9, p. 1].
  
In 2, we give the basic definitions, notations and certain standard results are obtained that indicate the unusual behavior of the algebra of all consequence operators defined on a set. In 4, some standard properties relative to subalgebras and chains in the set of all consequence operators are investigated. Finally, the entire last section is devoted to the foundations of the theory of nonstandard consequence operators defined on a nonstandard language.\par\bigskip

\indent {\bf 2. Basic concepts}\parm
Our notations and definitions for the standard theory of consequence operators are taken from references [2][6][12][13]$,$ and we now recall the most pertinent of these. Let L be any  nonempty set that is often called a {\it language}$,$ $\power {\r L}$ denote the power set of L and for any set X let $F(\r X)$ denote the finite power set of X (i.e. the set of all finite subsets of X.)\parm
DEFINITION 2.1 A mapping $\r C\colon \power{\r L} \to \power {\r L}$ is a consequence operator (or closure operator) if for each $\r X,\ \r Y \in 
\power {\r L}$\pars
\indent\indent (i) $\rm X \subset C(X) = C(C(X)) \subset L$ and if\pars
\indent\indent (ii) $\rm X \subset Y$$,$ then $\rm C(X) \subset C(Y).$\pars
\noindent A consequence operator C defined on L is said to be {\it finite} ({\it finitary}$,$ or {\it algebraic}) if it satisfies\pars
\indent\indent (iii) $\rm C(X) = \cup\{C(A)\mid A \in F(\r X)\}.$\par\medskip
REMARK 2.2 The above axioms (i) (ii) (iii) are not independent. Indeed$,$ 
(i)(iii) imply (ii).\pars
Throughout this entire article the symbol ``C'' with or without subscripts 
or with or without symbols juxtapositioned to the right will always denote a consequence operator. The only other symbols that will denote  consequence operators are ``I'' and ``U''. The symbol $\cal C$ [resp. ${\cal C}_f$] denotes the set of all consequence operators [resp. finite consequence operators] defined on $\power {\r L}.$\parm
DEFINITION 2.3. (i) Let\ \r I \ denote the identity map defined
on $\power {\r L}.$ \par
\noindent (ii) Let $\r U\colon \power {\r L} \to \power {\r L}$ be defined as follows: for each $\r X \in \power {\r L},\ \r U(\r X) = \r L.$\par
\noindent (iii) For each $\r C_1,\ \r C_2 \in {\cal C},$ define $\r C_1 \leq
\r C_2$ iff $\r C_1(\r X) \subset \r C_2(\r X)$ for each $\r X \in 
\power {\r L}.$
(Note that $\leq$ is obviously a partial order defined on $\cal C$.)\par  
\noindent (iv) For each $\r C_1,\ \r C_2 \in {\cal C}$$,$ define $\r C_1 \lor \r C_2
\colon \power {\r L} \to \power {\r L}$ as follows: for each $\r X \in 
\power {\r L},\ (\r C_1 \lor \r C_2)(\r X)= \r C_2(\r X)\cup \r C_2(\r X).$\par
\noindent (v) For each $\r C_1,\ \r C_2 \in {\cal C},$ define
$\r C_1 \land \r C_2 \colon \power {\r L} \to \power {\r L}$ as follows: for each $\r X \in \power {\r L},\ (\r C_2 \land \r C_2)(\r X) = \r C_1(\r X) \cap \r C_2 (\r X).$ \par
\noindent (vi) For each $\r C_1,\ \r C_2 \in {\cal C}$ define $\r C_2 \lor_w \r C_2\colon \power {\r L} \to \power {\r L}$ as follows: for each 
$\r X \in \power {\r L},\ (\r C_1 \lor_w \r C_2)(\r X) = \cap\{\r Y\mid
\r X \subset \r Y\subset \r L$ and $\r Y = \r C_1(\r Y)= \r C_2(\r Y)\}.$\pars 
Prior to defining certain special consequence operators notice that $\r I,\  \r U \in {\cal C}_f$ and that \ I [resp. U] \ is a lower [resp. upper] unit for the algebras $\langle {\cal C}, \leq \rangle$ and $\langle {\cal C}_f,\leq \rangle.$\parm

DEFINITION 2.4. Consider any  $\r X,\ \r Y \in \power {\r L}$. \par
\noindent (i) Define $\r C(\r X, \r Y)\colon \power {\r L} \to \power {\r L}$ as follows: let $\r A \in \power {\r L}.$ If $\r A \cap \r Y \not= \emptyset,$ then $\r C(\r X, \r Y)(\r A) = \r A \cup \r X.$ If $\r A \cap \r Y = \emptyset,$ then $\r C(\r X, \r Y)(\r A) = \r A.$\par
\noindent(ii) Define  $\r C'(\r X, \r Y) \colon \power {\r L} \to \power {\r L}$ as follows: let $\r A \in \power {\r L}.$ If $\r Y \subset \r A,$ then $\r C'(\r X, \r Y)(\r A) = \r A \cup \r X.$ If $\r Y \not\subset \r A,$ then $\r C'(\r X,\r Y)(\r A)= \r A.$\parm
{\bf THEOREM 2.5.} {\it For each $\r X,\ \r Y \in \power {\r L},\ \r C(\r X, \r Y) \in {\cal C}_f$ and $\r C'(\r X,\r Y) \in \cal C.$ If $\r Y\in F(\r L),$ then $\r C'(\r X,\r Y) \in {\cal C}_f.$}\pars 
PROOF. Let $\r X,\ \r Y,\ \r A \in \power {\r L}$ and consider $\r C(\r X,\r Y).$ If $\r A \cap \r Y \not= \emptyset,$ then $\r C(\r X,\r Y)(\r A)= \r A \cup \r X \supset \r A.$ If $\r A \cap \r Y = \emptyset,$ then $\r C(\r X,\r Y)(\r A)= \r A.$ Hence, for each $\r A \in \power {\r L},\ \r A \subset \r C(\r X,\r Y)(\r A).$ Assume $\r A \cap \r Y \not= \emptyset.$ Then $\r C(\r X,\r Y)(\r C(\r X,\r Y)(\r A))= \r C(\r X,\r Y)(\r A\cup \r X) = \r A \cup \r X = \r C(\r X,\r Y)(\r A),$ since $(\r A \cup \r X)\cap \r Y \not= \emptyset.$ If $\r A \cap \r Y =\emptyset,$ then $\r C(\r X,\r Y)(\r C(\r X,\r Y)(\r A)) = \r C(\r X,\r Y)(\r A).$ Thus, axiom (i) of definition 2.1 holds. Let $\r A \subset \r H \subset \r L.$ If $\r A \cap \r Y \not=\emptyset,$ then $\r H\cap \r Y\not= \emptyset$ implies that $\r C(\r X, \r Y)(\r A) = \r A \cup \r X \subset \r H \cup \r X = \r C(\r X,\r Y)(\r H).$ Assume that $\r A\cap \r Y = \emptyset.$ Then $\r C(\r X, \r Y)(\r A) = \r A.$ If $\r H \cap \r Y \not=\emptyset,$ then $\r A \subset \r H \cup \r X = \r C(\r X,\r Y)(\r H).$ If $\r H\cap \r Y = \emptyset,$ then $\r A \subset \r H = \r C(\r X,\r Y)(\r H).$ Thus, in all cases, $\r C(\r X,\r Y)(\r A) \subset \r C(\r X,\r Y)(\r H)$ and axiom (ii) holds. Let $\r A \cap \r Y \not=\emptyset$ and $ \r x \in \r C(\r X,\r Y)(\r A)= \r A \cup \r X.$ If $ \r x \in \r A,$ then $\r C(\r X, \r Y) (\{\r x\}) = \{\r x\} \cup \r X$ or $\{ \r x\}.$ Hence, in this case, $\r x \in \r C (\r X,\r Y)(\{\r x\}).$ Suppose that $\r x \in \r X.$ Then there exists some $\r y \in \r A \cap \r Y$ and $\r x \in \r C(\r X, \r Y)(\{\r y\}) = \{\r y\}\cup \r X \subset \r A \cup \r X.$ Consequently, if $\r A \cap \r Y \not=\emptyset$ and $\r x \in \r C(\r X, \r Y)(\r A)$, then there is some $\r F \in F(\r A)$ such that $\rm x \in C(X,Y)(F).$ Consider the case where $\r A \cap \r Y = \emptyset.$ If $\rm A = \emptyset,\ C(X,Y)(A) =\emptyset =  \bigcup\{C(X,Y)(F)\mid F \in {\it F}(\emptyset)\}.$ Let $\rm A \not=\emptyset,$ then $\r x \in \r C(\r X, \r Y)(\r A)=\r A$ implies that  $\r x \in \r C(\r X,\r Y)(\{\r x\})$ and $\{\r x\} \in F(\r A).$ Hence, in general, if $\r x\in \r C(\r X, \r Y)(\r A),$ then $ x \in \bigcup \{\r C(\r X, \r Y)(\r F)\mid \r F \in F(\r A)\}.$ Since $\r C(\r X,\r Y)(\r F) \subset \r C(\r X,\r Y)(\r A)$ for each $\r F \in F(\r A),$ then  it follows that (iii) holds.\pars
Consider $\r C'(\r X,\r Y),$ let $\r A \in \power {\r L}$ and assume that $\r Y \subset \r A.$ Then $\r C'(\r X, \r Y)(\r A)= \r A \cup \r X \supset \r A.$ Moreover, $\r C'(\r X, \r Y)(\r C'(\r X, \r Y)(\r A)) = \r C'(\r X, \r Y)(\r A \cup \r X)= \r A \cup \r X = \r C'(\r X, \r Y)(\r A).$ If $\r Y \not\subset \r A,$ then $\r C'(\r X, \r Y)(\r A) = \r A$ and $\r C'(\r X, \r Y)(\r C'(\r X, \r Y)(\r A))= \r C'(\r X, \r Y)(\r A).$ Thus axiom (i) holds. The fact that if $\r A \subset \r H \subset \r L,$ then $\r C'(\r X, \r Y)(\r A) \subset \r C'(\r X, \r Y)(\r H)$ follows easily and (ii) holds.\pars
Assume that $\r Y \in F(\r L),\ \r Y \subset \r A$ and $\r x \in \r C'(\r X, \r Y)(\r A) = \r A \cup \r X.$ If $\r x \in \r X,$ then $\r x\in \r C'(\r X, \r Y)(\r Y) = \r Y \cup \r X \subset \r A \cup \r X.$  If $ \r x \in \r A,$
then $\r x \in \r C'(\r X, \r Y)(\r Y \cup \{\r x\}) = \r Y \cup \{\r x\} \cup \r X \subset \r A \cup \r X.$ But $\r Y \cup \{\r x\} \in F(\r A).$ Hence, in this case, $\r x \in \r C'(\r X,\r Y)(\r F),$ where $\r F \in F(\r A).$ Finally, let $\r Y \not\subset \r A.$ If $\rm A = \emptyset,\ C'(X,Y)(A) = \emptyset = \bigcup\{C'(X,Y)(F)\mid F \in {\it F}(\emptyset)\}.$ Assume $\rm A \not= \emptyset$ and $ \r x \in \r C'(\r X, \r Y)(\r A) = \r A.$ Then $\r x \in \r A.$ If $\r A \subset \r Y,\ \r A \not= \r Y,$ then $\r x \in \r C'(\r X, \r Y)(\{\r x\})= \{\r x\}.$ Otherwise,
$\r A \not\subset \r Y$ and there exists some $\r z \in \r A$ such that $\r z \not\in \r Y.$ In which case, $\r Y \not\subset \{\r x, \r z\}$ and, hence, $\r x \in \r C'(\r X, \r Y)(\{\r x, \r z\})= \{ \r x, \r z\} \in F(\r A).$ Therefore,  $\r C(\r X, \r Y)(\r A) \subset \bigcup\{\r C'(\r X,\r Y)(\r F)\mid \r F \in F(\r A)\}.$ This result and axiom (ii) imply that axiom (iii) holds. This completes the proof.\pars 
Recall that $\r C \in \cal C$ is {\it axiomless} if $\r C(\emptyset) = \emptyset$ and axiomatic otherwise [8]. Note that for any $\r X,\ \r Y \in \power {\r L},\ \r C(\r X,\r Y)$ is axiomless, and if $\r X = \emptyset$ or $\r Y \not=\emptyset,$ then $\r C'(\r X,\r Y)$ is axiomless.\parm
{\bf LEMMA 	2.6.} {\it Let $\r C \in \cal C$ be axiomatic. Then there exists some $\r x \in \r L$ such that $\r C(\r L - \{\r x\}) = \r L.$}\pars
PROOF. Assume that there does not exist some $\r x \in \r L$ such that $\r C(\r L - \{\r x\}) = \r L.$ Then for each $\r y \in \r L,\ \r L - \{\r y \} \subset \r C(\r L - \{\r y\}) \subset \r L$ implies that $\r L - \{\r y\} = \r C(\r L - \{\r y\})$ from axiom (i). But axiom (ii) yields that $\r C(\emptyset) = \r C(\bigcap\{\r L - \{ \r y\}\mid \r y \in \r L\}) \subset \bigcap \{\r C(\r L - \{ \r y \})\mid \r y \in \r L \} = \bigcap\{\r L - \{ \r y \}\mid \r y \in \r L \} = \emptyset.$ Thus $\r C$ would be axiomless and this contradiction completes the proof. \pars
Recall that a member $\r C_1$ in the algebra $\langle {\cal C}, \leq \rangle$ {\it covers} $\r C_2 \in {\cal C},$ if $\r C_2 < \r C_1$ and 
there does not exist some $\r C_3 \in \cal C$ such that $\r C_2 < \r C_3 < \r C_1.$ A set ${\cal B} \subset \cal C$ {\it densely} covers ${\cal E} \subset \cal C$ if for each $\r C \in \cal B$ there exists some $\r C_1 \in \cal E$ such that $\r C_1 \leq \r C.$ Recall that $\r C \in \cal C$ is an {\it atom} if $\r C $ covers $\r I$ and $\langle {\cal C},\leq, \r I \rangle$ is {\it atomic} with ${\cal E}\subset \cal C$ the set of atoms if each member of $\cal E$ covers \ I\ and for each $\r I\not=\r C \in {\cal C}$ there is some $\r C_1 \in \cal E$ such that $\r C_1 \leq \r C.$ Let 
${\cal E}_0 = \{\r C'(\{\r x\},\r L- \{ \r x\})\mid \r x \in \r L\}.$ Notice that $\r I \notin {\cal E}_0,$ since for $\r x \in \r L,\ \r C'(\{\r x\},\r L - \{\r x \})(\r L - \{ \r x\}) = \r L,$  and that each member of ${\cal E}_0$ is axiomless if $\r L$ has more than one member. The next result shows that $\langle {\cal C}, \leq \rangle$ is almost atomic.\parm
{\bf THEOREM 2.7.} {\it For $\langle {\cal C}, \leq \rangle,$ the set all axiomatic consequence operators ${\cal C}_A$ covers ${\cal E}_0$ and each member of ${\cal E}_0$ is a atom.}\pars
PROOF. We first show that each member of ${\cal E}_0$ is an atom. Let $\r x\in \r L$ and assume that there exists some $\r C_1 \in \cal C$ such that $\r C_1 < \r C'(\{ \r x\}, \r L - \{ \r x \}) = \r C'.$ Assume that $\r B \in \power {\r L},\ \r C_1(\r B) \subset \r C'(\r B)$ and  $\r C_1(\r B) \not= \r C'(\r B).$ Suppose that $\r L - \{ \r x \} \subset \r B.$ Then $\r C'(\r B) = \r L$ and $\r L - \{ \r x \} \subset \r B \subset \r C_1(\r B) \not= \r L.$  Hence, $\r L - \{ \r x \} = \r B = \r C_1(\r B).$ Now suppose that $\r L - \{ \r x \} \not\subset \r B \not= \r L.$ Then $\r B \subset \r C_1 (\r B) \subset \r C'(\r B) = \r B.$ However, this contradicts $\r C_1(\r B) \not= \r C'(\r B).$ Thus for any $\r B \in \power {\r L}$ such that $\r C_1 (\r B) \not= \r C'(\r B),$ it follows that $\r L -\{\r x\} \subset \r B$ and $\r C_1(\r B) = \r B$ and there is only one such set with these properties, the set is $\r L - \{\r x\}$ since $\r C_1(\r L) = \r C'(\r L) = \r L.$ Therefore, it must follow that $\r L - \{ \r x \} = \r B = \r C_1(\r B)= \r I(\r B)$ and $\r C'(\r L) = \r L$ in order that $\r C_1(\r B) \subset \r C'(\r B)$ and $\r C_1(\r B) \not= \r C'(\r B).$ Further, in general, $\r I (\r B) \subset \r C'(\r B) = \r L$ and $\r I(\r B) \not= \r C'(\r B).$ Now if  $ \r A \in \power {\r L}$ and $ \r B \not=\r A \not= \r L,$ then $\r L - \{ \r x \}\not\subset \r A$ implies that $\r C_1(\r A) = \r C'(\r A)= A.$  Finally, if $\r A = \r L,$ then $\r C_1(\r L) = \r C'(\r L) = \r L = \r I (\r L).$ Consequently, $\r C_1 = \r I,$ $\r C_1 < \r C'$ and there is no $\r C \in \cal C$ such that $\r C_1 < \r C < \r C'.$ Hence, $C'$ is an atom.\pars 
We now easily show that ${\cal C}_A$ densely covers ${\cal E}_0.$ Let $\r C \in {\cal C}_A.$ Then from Lemma 2.6, there exists some $\r x \in \r L$ such that $\r C(\r L - \{ \r x \}) = \r L.$ Then letting $\r C'(\{\r x\},\r L - \{ \r x\}) = \r C',$ it follows that $\r C'(\r L -\{\r x\}) = \r L = \r C(\r L -\{ \r x \}).$ If $\r A = \r L,$ then $\r C'(\r A) = \r C(\r A)$ and if $ \r A \subset \r L,\ \r A \not= \r L - \{ \r x\},\  \r A \not= \r L,$ then $\r C'(\r A) = \r A \subset \r C(\r A).$ Hence, $\r C' \leq \r C.$ This completes the proof.\pars
It is not difficult to show that $\langle {\cal C},\leq \rangle$ [13] and $\langle {\cal C}_f, \leq \rangle$ are both closed under $\land$,
where $\land$ is (v) of definition 2.3, which, obviously, would be the {\it meet} operation associated with $\leq.$ However, as will be shown by a simple example, it is rare that these algebras are closed under $\lor,$ which if either is so closed, then $\lor$ would be the {\it join} operation. W\'ojcicki was the first to recognize that not only are these algebras closed under $\lor_w,$ and $\lor_w$ is the {\it join} operation, but $\langle {\cal C},\land,\lor_w, \r I ,\r U\rangle$ is also complete [13, p. 276]. Unfortunately, $\langle {\cal C}_f,\land,\lor_w, \r I ,\r U\rangle$  is not complete [2, p. 180]. For a simple proof that $\langle {\cal C},\land,\r I, \r U\rangle$ is meet-complete (and thus complete) see [13, p. 276]. Using the fact that $\langle {\cal C},\leq \rangle$ is a meet semi-lattice, it follows easily that $\langle {\cal C}_f, \leq \rangle$ is also a meet semi-lattice. We need only show that, for $\langle {\cal C}_f, \leq \rangle,$ $\land$ satisfies axiom (iii). Let $\r C_1,\ \r C_2 \in {\cal C}_f.$ Then for each $\r A \in \power {\r L},\ (\r C_1 \land \r C_2)(\r A) = \r C_1(\r A) \cap \r C_2 (\r A)= (\cup \{\r C_1 (\r X)\mid 
\r X \in F(\r A)\})\bigcap (\cup \{\r C_2 (\r X)\mid \r X \in F(\r A)\})=\bigcup \{\r C_1(\r X) \cap \r C_2(\r X)\mid \r X \in F(\r A)\}=\bigcup \{(\r C_1\land \r C_2)(\r X)\mid \r X \in F(\r A)\}$ and (iii) holds.\parm
EXAMPLE 2.8. Certain subsets of  $\langle {\cal C},\leq \rangle$ and $\langle {\cal C}_f, \leq \rangle$ may be closed under $\lor,$ but, in general, there are members such that the $\lor$ operator does not yield a consequence operator.  Let $\r L$ have 3 or more members. Define $\r S \colon \power {\r L} \to \power {\r L}$ by letting $\emptyset \not= \r M \subset \r L,\ \vert \r L - \r M \vert \geq 2.$ Let $\r S(\emptyset ) = \r M$ and $\r b \in \r L - \r M.$ For each $\r A \in \power {\r L},$ if $\r b \in \r A,$ then let $\r S(\r A) = \r L$; if $\r b \notin \r A,$ let $\r S(\r A) = \r M \cup \r A.$ It follows easily that $\r S \in {\cal C}_A \cap {\cal C}_f.$ Consider $\r C'(\{ \r b\},\emptyset) \in {\cal C}_f.$ Then $(\r C'(\{ \r b\},\emptyset) \lor \r S)(\r C'(\{ \r b\},\emptyset) \lor \r S)(\emptyset) = (\r C'(\{ \r b\},\emptyset) \lor \r S)(\r C'(\{\r b\},\emptyset)(\emptyset) \cup \r S(\emptyset)) = (\r C'(\{ \r b\},\emptyset) \lor \r S)(\{\r b\} \cup \r M) = (\r C'(\{ \r b \}, \emptyset)(\{ \r b \} \cup \r M) \cup \r S(\{ \r b \} \cup \r M)= \{ \r b\} \cup \r M \cup \r L = \r L \not= (\r C'(\{\r b\},\emptyset)\lor \r S)(\emptyset) = \{ \r b\} \cup \r M.$ Thus $\r C'(\{ \r b\},\emptyset) \lor \r S$ is not a consequence operator. \pars
Observe that if $\r L$ is a standard formal propositional or predicate language and $\r S'$  the propositional or predicate consequence operator respectively, then even though $\r S'$ is not the same operator as defined in example 2.8 the presence of the formula $\r b = \r P \land (\neg \r P)$ 
will also yield that $\r C'(\{\r b\}, \emptyset)) \lor \r S' \notin {\cal C}.$ 
Simply substitute $\r S'$ for $\r S$ with this formula as the $\r b$ and notice that $\r L - \r S'(\emptyset)$ is denumerable, $\r b \in (\r L - \r S'(\emptyset)),$ and $\r S'(\{\r b \})\cup \r S'(\emptyset) = L.$ \par\bigskip
{\bf 3. Subalgebras}\parm 
Subalgebras of $\langle {\cal C},\land,\lor_w, \r I, \r U \rangle$ and  
$\langle {\cal C}_f,\land,\lor_w, \r I, \r U \rangle$ have been studied to a certain extent and appear to be the most appropriate area for further investigation. It is known that there are sublattices of $\langle {\cal C},\land,\lor_w, \r I, \r U \rangle$ that are atomic and coatomic [2, p. 179]. We first show that there are complete and distributive sublattices of 
$\langle {\cal C}_f,\land,\lor_w, \r I, \r U \rangle,$ where $\lor = \lor_w.$ Moreover, such sublattices need not be atomic.\parm
{\bf THEOREM 3.1.} {\it For each $\r B \in \power {\r L},$ let ${\cal C}(\r B) = \{\r C(\r X,\r B)\mid \r X \in \power {\r L} \}.$ Then $\langle {\cal C}(\r B), \land,\lor,\r I , \r C(\r L,\r B)\rangle$ is a complete and distributive sublattice of $\langle {\cal C}_f,\land,\lor_w, \r I, \r U \rangle$. If there exists nonempty $\r A,\ \r B\in \power {\r L},$ such that $\r A \not= \r L, \  \r B \subset \r A,$ and $\r B \not= \r A,$ then $\langle {\cal C}(\r B), \land,\lor,\r I , \r C(\r L,\r B)\rangle$ is not a chain.}  \pars
PROOF. Let $\r A,\ \r B,\ \r X,\ \r Y \in \power {\r L}.$ Note that $\r C(\emptyset,\r B) = \r I.$ Let $\cal H$ be any nonempty subset of ${\cal C}(\r B).$ Then there exists some nonempty ${\cal A} \subset \power {\r L}$ such that ${\cal H}= \{\r C(\r A,\r B)\mid \r A \in {\cal A}\}.$ We first show that $\inf ({\cal H}) = \r C( \cap {\cal A}, \r B).$ If $\r B \cap \r X \not= \emptyset,$ then $\r C(\cap {\cal A},\r B) (\r X)= (\cap {\cal A}) \cup \r X \subset \r C(\r A, \r B)(\r X) = \r A \cup \r X$ for each $\r A \in {\cal A}.$ If $\r B \cap \r X = \emptyset,$ then $\r C(\cap {\cal A},\r B)(\r X) = X = \r C(\r A, \r B)(\r X)$ for each $\r A \in {\cal A}.$ Thus $\r C(\cap {\cal A}, \r B) \leq \r C(\r A, \r B)$ for each $\r A \in {\cal A}.$ Hence, $\r C(\cap {\cal A}, \r B)$ is a lower bound for ${\cal H}.$ Let $\r C(\r Y,\r B)$ be a lower bound for $\cal H$. If $\r B \cap \r X \not=\emptyset,$ then $\r Y \cup \r X \subset \r A \cup \r X$ for all $\r A \in \cal A$ yields that $\r Y \cup \r X \subset \bigcap\{(\r A \cup \r X)\mid \r A \in {\cal A}\} = (\cap {\cal A}) \cup \r X = \r C(\cap  {\cal A},\r B)(\r X).$ If $\r B \cap \r X = \emptyset,$ then $\r C(\r Y,\r B)(\r X) = \r X = \r C(\r A,\r B)(\r X)= \r C(\cap {\cal A},\r B)(\r X)$ where $\r A \in {\cal A}.$ Thus $\r C(\r Y, \r B) \leq \r C(\cap {\cal A},\r B)$
implies that $\inf({\cal H}) = \r C(\cap {\cal A}, \r B).$ We next show that $\sup({\cal H}) = \r C(\cup {\cal A},\r B).$ But, first, we show that $\lor = \lor_w.$ Let $\r A_1 \in \power {\r L}$ and assume that $\r B \cap \r X = \emptyset.$ Then $(\r C(\r A,\r B) \lor \r C(\r A_1,\r B))(\r X) = \r X.$ Now from the definition of $\lor_w$ [13, p. 176], $(\r C(\r A, \r B) \lor_w \r C(\r A_1, \r B))(\r X) = \bigcap \{ \r Y \mid \r X \subset \r Y = \r C(\r A, \r B)(\r Y) = \r C(\r A_1, \r B) (\r Y)\} = X$ since $\r C(\r A, \r B)(\r X)=\r C(\r A_1,\r B)(\r X) = \r X$ in this case. Now let $\r B \cap \r X \not= \emptyset.$ Then $(\r C(\r A, \r B)\lor \r C(\r A_1, \r B))(\r X) = \r A \cup \r A_1 \cup \r X.$ If $\r X \subset \r C(\r A, \r B)(\r Y) = \r Y= \r C(\r A_1,\r B)(\r Y),$ then $\r Y \cap \r B \not=\emptyset$ yields that $\r A \cup \r Y = Y = \r A_1 \cup \r Y.$ Hence, letting $\r Y_1 = \r A \cup \r A_1 \cup \r X \subset \r Y,$ then $ \r C(\r A, \r B)(\r Y_1) = \r Y_1 = \r C(\r A_1, \r B)(\r Y_1).$ Thus $(\r C(\r A, \r B) \lor_w \r C(\r A_1, \r B))(\r X) = \bigcap \{ \r Y \mid \r X \subset \r Y = \r C(\r A, \r B)(\r Y)= \r C(\r A_1, \r B)(\r Y)\} = \r A \cup \r A_1 \cup \r X.$ Consequently, 
$\lor = \lor_w.$ It now follows, in like manner, that $\sup({\cal H}) = \r C(\cup {\cal A}, \r B).$ The fact that this complete  sublattice is distributive follows from the fact that $\cup$ and $\cap$ are distributive.\pars
For the final part, assume that nonempty $\r A, \ \r B\in \power {\r L},\ \r  A \not= \r L,\ \r B \subset \r A$ and $\r B \not= \r A.$ There is some $\r x \in \r L - A.$ Let $\r D = \{\r x\}.$ Then $\r C(\r A, \r B)(\r B) = \r A \cup \r B = \r A$ and $\r C(\r D, \r B)(\r B) = \r B \cup \r D.$ However, $\r A \not\subset \r B \cup \r D$ and $\r B \cup \r D \not\subset \r A$ imply that $\r C(\r A, \r B)$ and $\r C(\r D, \r B)$ are not comparable. Thus in this case, $\langle {\cal C}(\r B), \land, \lor, \r I, \r C(\r L, \r B) \rangle$ is not a chain and the proof is complete. 
\pars
Using the collection of axiomless consequence operators defined in theorem 3.1, we show, in general, that the algebras $\cal C$ and ${\cal C}_f$ do not have the descending chain condition.\parm
EXAMPLE 3.2. For infinite $\r L,$ we show that for a specific $\r B$, $\langle {\cal C}(\r B), \land,\lor,\r I , \r C(\r L,\r B)\rangle$ contains a $\langle {\cal C}_f,\leq\rangle$ chain that does not satisfy the descending chain condition. Hence, $\langle {\cal C}_f,\leq\rangle$ and $\langle {\cal C},\leq\rangle$ would not satisfy the descending chain condition. There exists an injection $\rm F\colon \nat \to L$. Let $\rm f(0) = x_0$ and $
\rm B =\{x_0\}.$ For each $\rm n \in \nat,\ n\geq 1,$ let $\rm F_n =L - f[[1,n]]$ and $\rm C_n = C(F_n,B).$ 
Notice that no $\rm F_n =\emptyset.$ Let $\rm {\cal H} = \{C_n\mid n \geq 1\}.$ Let distinct $\rm C_k,\ C_m \in {\cal H}.$  Then either $\rm k < m$ or $\rm m < k$. Assume that $\rm k < m.$ In general, let $\rm C_n \in {\cal H},\ X \in \power {L}.$ If $\rm x_0 \in X,\ C_n(X) = F_n \cup X.$ If $\rm x_0 \notin X,\ C_n(X) = X.$ Thus, for any $\rm n < m, \ \emptyset\not= F_m {\lower4pt\hbox{$\buildrel \subset \over{\not=}$}} F_n$ implies that $\rm C_m {\lower4pt\hbox{$\buildrel < \over{\not=}$}} C_k.$ In like manner, for $
\rm m < k.$ Hence, any two members of $\cal H$ are comparable with respect to $\langle {\cal C}_f,\leq\rangle.$  Also for each $\rm C_n,\ \rm C_n(\{x_0,x_1\}) = F_n \cup \{x_0\} = F_n\not= \{x_0,x_1\}.$ Therefore, no $\rm C_n = I.$ This yields that for $\rm n\in \nat,\ n\geq 1,\ I \not= C_{n+1} {\lower4pt\hbox{$\buildrel < \over{\not=}$}} C_n$ and this completes the example. Notice that by taking $\rm A =\{x_0,x_1\},\ \rm {\cal C}(B)$ is not a chain.\pars 
Since $\langle {\cal C},\land,\lor_w,\r I, \r U\rangle$ contains distributive and complete sublattices that are not chains, where $\lor = \lor_w,$ then a natural question to ask is whether or not such sublattices can have any other Boolean type structures?\parm
{\bf THEOREM 3.3.} {\it Let $\r C,
 \r C_1 \in {\cal C}$ and 
 $\r I\ {\lower4pt\hbox{$\buildrel < \over{\not=}$}}\ \r C \ {\lower4pt\hbox{$\buildrel < \over{\not=}$}}\ \r C_1.$ Then $\r C' \in \cal C$ is a complement relative to $\r C$ and with respect to $\lor$ iff for each $A \in \power {L},\ \r C'(\r A) = (\r C_1(\r A) - \r C(\r A)) \cup \r A.$}\pars
PROOF.  Let $\r C,\ \r C_1 \in {\cal C},$ and $\r I\ {\lower4pt\hbox{$\buildrel < \over{\not=}$}}\ \r C \ {\lower4pt\hbox{$\buildrel < \over{\not=}$}}\ \r C_1.$ The $\r C'$ is a relative complement for $\r C$ iff $\r C \lor \r C' = \r C_1$ and $\r C \land \r C' = \r I$ iff for each $\r A \in \power {\r L},$ (1) $\r C(\r A) \cup \r C'(\r A) = \r C_1 (\r A)$ and (2)
$\r C(\r A) \cap \r C'(\r A) = \r A.$ For the necessity, assume (1) and (2) and let $\r y \in  (\r C_1(\r A) - \r C(\r A)) \cup \r A.$ If $\r y \in \r A,$ then $\r y \in \r C'(\r A)$ by (2). If $\r y \notin \r A,$ then $\r y \in \r C_1(\r A) - \r C(\r A)$ implies that $\r y \in \r C_1(\r A)$ and $\r y \notin \r C(\r A)$; which further implies that $\r y \in \r C'(\r A)$ from (1). One the other hand, if $\r y \in \r C'(\r A)$ and $\r y \not \in \r A,$ then by (1) and (2), $\r y \in \r C_1 (\r A) - \r C(\r A)$ implies that 
$\r y \in  (\r C_1(\r A) - \r C(\r A)) \cup \r A.$ Obviously, if $\r y \in A,$ then $\r y \in (\r C_1(\r A) - \r C(\r A)) \cup \r A.$ Thus if (1) and (2) holds, then $ \r C'(\r A) =  (\r C_1(\r A) - \r C(\r A)) \cup \r A$ for each $\r A\in \power {\r L}$\pars
For the sufficiency, let $\r C'(\r A) = (\r C_1(\r A) - \r C(\r A)) \cup \r A$ for each $\r A \in \power {\r L}.$ Then $\r C(\r A) \cup \r C'(\r A) = \r C_1(\r A) \cup \r A = \r C_1(\r A)$ and $ \r C'(\r A) \cap \r C(\r A) =((\r C_1(\r A) - \r C(\r A)) \cup \r A)\cap \r C(\r A) = \emptyset\cup \r A = \r A$ and this completes the proof.\pars 
EXAMPLE 3.4. We show that $\cal C$ is not closed under composition of maps even though such a composition is always defined. Let $\r S$ be the consequence operator defined in example 2.8. Then $\emptyset \not= \r S(\emptyset) = \r M \not= \r L,$ $\r b \in \r L - \r M$ and there exists some $\r a \in \r L - \r M$ such that $\r a \not= \r b.$ Hence, denoting composition by juxtaposition, it follows that $(\r C'(\{ \r b\},\emptyset)\r S)\colon \power {\r L} \to \power {\r L}$ and $(\r C'(\{ \r b\}, \emptyset)\r S)(\r M) = \r C'(\{\r b\},\emptyset)(\r M) = \r M \cup \{\r b\}.$ However, $(\r C'(\{ \r b\},\emptyset)\r S)(\r M\cup \{ \r b\}) = \r C'(\{\r b\},\emptyset)(\r L) = \r L.$ Therefore, $(\r C'(\{ \r b\},\emptyset)\r S)(\r C'(\{ \r b\},\emptyset)\r S)(\r M) = \r L.$ It follows that this composition is not a consequence operator. \pars
Even though example 3.4 shows that $\cal C$ is not closed under composition  the next proposition gives a strong characterization for chains in terms of composition.\parm
{\bf THEOREM 3.5.} {\it Let ${\cal A} \subset \cal C.$ Then $\cal A$ is a chain in $\langle {\cal C}, \leq \rangle$ iff for each $\r C, \r C' \in \cal A$ either the composition $\r C'\r C = \r C'$ or $\r C\r C' = \r C.$}\pars
PROOF. Let $\r C,\ \r C' \in \cal A$ and assume that $\cal A$ is a chain in 
$\langle {\cal C}, \leq \rangle.$ Suppose that $\r C \leq \r C'.$ Then for each $\r B \in \power {\r L},\ \r B \subset \r C(\r B)\subset \r C'(\r B)$ implies that $\r C'(\r B) \subset \r C'\r C(\r B) \subset \r C'(\r C'(\r B)) = \r C'(\r B).$ Hence, $\r C'\r C= \r C'.$ In like manner, if $\r C' \leq \r C,$ then $\r C\r C' = \r C.$\pars
Conversely, let $\r C,\ \r C' \in \cal A$ and $ \r C'\r C = \r C'.$ Then for each $\r B\ \in \power {\r L},\ \r C(\r B)\subset \r C'(\r C(\r B)) = (\r C'\r C)(\r B) = \r C'(\r B).$ Thus $\r C \leq \r C'$. In like manner, if $\r C\r C' = \r C,$ then $\r C' \leq \r C$ and this completes the proof. \pars
In the next section, our attention is often restricted to chains in $\langle  {\cal C}_f,\leq \rangle.$ We first embed $\langle {\cal C},\leq\rangle$ into a non-standard structure and investigate nonstandard bounds for various chains.\parm  
  
{\bf 4. Nonstandard Consequence Operators}\parm
Let ${\cal A}$ be a nonempty finite set of symbols. It is often convenient to assume that $\cal A$ contains a symbol that represents a blank space. As usual any nonempty finite string of symbols from $\cal A,$ with repetitions$,$ is called a {\it word} [10$,$ p.222]. A word is also said to be an (intuitive) {\it readable sentence} [5$,$ p. 1]. We let \ W\ be the intuitive set of all words created from the {\it alphabet} $\cal A$. Note that in distinction to the usual approach$,$ \ W \ does not 
contain a symbol for the empty word. \pars
We accept the concept delineated by Markov [4]$,$ the so-called  ``abstraction 
of identity$,$'' and say that $\r w_1,\ \r w_2 \in \r W$ are ``equal'' if they 
are composed of the same symbols written in the same intuitive order (left to 
right). The {\it join} or juxtaposition operation between $\r w_1,\ \r w_2 \in 
\r W$ is the concept that yields the string $\r w_1 \r w_2$ or $\r w_2 \r 
w_1$. Thus \ W\ is closed under join. Notice that we may consider a 
denumerable formal language as a subset of \ W. (By adjoining a new symbol not in $\cal A$ and defining it as the unit, $\r W$ becomes a free monoid generated by the set ${\cal A} \cup \{\rm new\ symbol\}.$)\pars
Since \ W\ is denumerable$,$ then there exists an injection
$i\colon \r W \to \nat.$ Obviously$,$ if we are working with a formal language that is a subset of \ W$,$ 
then we may require $i$ restricted to a formal language to be a G\"odel 
numbering. Due to the join operation$,$ a fixed member of \ W\ that contains two 
or more 
distinct symbols can be represented by various {\it subwords} that are joined together to yield the given fixed word. The word ``mathematics'' is generated by the join of $\r w_1 = {\rm math},\ \r w_2 = \r e,\ \r w_4 = {\rm mat},\ \r w_4 = {\rm 
ics}.$ This word can also be formed by joining together 11 not necessarily distinct members of \ W. \pars
Let $i[\r W] = T$ and for each $n \in \nat,$ let $T^n = T^{[0,n]}$ denote the set of all mappings from $[0,n]$ into $T.$ Each element of $T^n$ is called a {\it partial sequence}$,$ even though this definition is a slight restriction of the usual one that appears in the literature.
Let $f \in T^n, n>0.$ Then the {\it order induced by f} is the simple inverse order determined by $f$ applied to the simple order on $[0,n].$ Formally$,$ for each $f(j),\ f(k)\in f[[0,n]],$ define $f(k)\leq_f f(j)$ iff $j\leq k,$ where $\leq$ is the simple order for $\nat$ restricted to $[0,n].$ In general$,$ we will not use this notation $\leq_f$ but rather we will indicate this (finite) order in the usual acceptable manner by 
writing the symbols
$f(n),f(n-1),\ldots,f(0)$ from left to right . Thus we symbolically let $f(n)\leq_f f(n-1)\leq_f \cdots \leq_f f(0)= f(n)f(n-1)\ldots f(0).$\pars
Let $f\in t^n.$ Define $\r w_f \in \r W$ as follows: $\r w_f = 
(i^{-1}(f(n)))(i^{-1}(f(n-1)))\cdots (i^{-1}(f(0))),$ where the operation indicated by juxtaposition is the join. We now define a relation on 
$P =\cup \{T^n\mid n \in \nat \}$ as follows: let $f,\ g\in P.$ Then for 
$f\in T^n$ and $g \in T^m,$ define $f \sim g$ iff $(i^{-1}(f(n)))\cdots (i^{-1}f(0)))=(i^{-1}(g(m)))\cdots (i^{-1}(g(0))).$ It is obvious that 
$\sim$ is an equivalence relation on $P.$ For  each $f \in P,$ $[f]$ denotes 
the equivalence class under $\sim$ that contains $f.$ Finally$,$ let
${\cal E}=\{[f]\mid f\in P\}.$ Observe that for each $[f] \in \cal E$ there exist $f_0,\ f_m \in [f]$ such that 
$f_0 \in T^0,\ f_m \in T^m$ and if there exists some $k\in \nat$ such that $0<k<m,$ then there exists some
$g_k \in [f]$ such that $g_k \in T^k$ and if $j\in \nat$ and $j >m,$ then there does not exist $g_j \in T^j$ such that $g_j \in [f].$ If we define the {\it size} of a word $\r w \in \r W$ (size(w)) to be the number of not necessarily distinct symbols counting left to right that appear in  \ W$,$ then the size(w) = $m+1.$ For each $\r w \in \r W,$ there is $f_0 \in T^0$ such that
$\r w = i^{-1}(f_0(0))$ and such an $f_m \in [f_0]$ such that size (w) = $m+1.$
On the other hand$,$ given $f \in P,$ then there is a $g_0 \in [f]$ such that 
$(i^{-1}(g_0(0))) \in \r W.$ Of course$,$ each $g \in [f]$ is interpreted to be 
the word $(i^{-1}(g(k)))\ldots (i^{-1}(g(0))).$\pars
Each $[f] \in \cal E$ is said to be a (formal) word or (formal) {\it readable
sentence.} All the intuitive concepts$,$ definitions and results relative to consequence operators defined for $\r A  \in \power {\r W}$ are now passed
to $\power {{\cal E}}$ by means the map $\theta(i(\r w))= [f_0]$.
In the usual manner$,$ the map $\theta$
is extended to subsets of each $\r A \in \power {\r W},$ n-ary relations and the like. For example$,$ let $\r w \in \r A \in \power {\r w}.$ Then there exists $f_w \in P$ such that $f_w \in T^0$ and $f_w(0)= i(\r w).$ Then 
$\theta(i(\r w)) = [f_w].$ In order to simplify notation$,$ the images of the extended 
  $(\theta\, i)$ composition will often be indicated by bold notation with the exception  
of customary relation symbols which will be understood relative to the context. 
For example$,$ if \ S \ is a subset of \ W$,$ then we write $\theta (i[\r S]) = \b S.$\pars

Let $\cal N$ be a superstructure constructed from the set $\r W\cup \nat$ as its set of atoms. (12/16/12 the ground set has now been so expanded so any symbols in W that are used for natural numbers are different from those in $\nat$.) Our standard structure is ${\cal M} = ({\cal N}, \in, = ).$ Let 
$\Hyper {\cal M}= (\Hyper {\cal N},\in,=)$ be a nonstandard and elementary extension of $\cal M$$.$ Further, $\Hyper {\cal M}$ is an enlargement.\pars

For an alphabet $\cal A$$,$ there exists $[g] \in \Hyper {\cal E}-{\cal E}$ such that there are only finitely many standard members of $\nat$ in the range of $g$ and these 
standard members injectively correspond to alphabet symbols in $\cal A$ [5$,$ p. 24]. On the other hand$,$ there exist $[g'] \in \Hyper {\cal E}-{\cal E}$ such that the range of $[g']$ does not correspond in this manner to elements in $\cal A$ [5$,$ p. 90]. \pars
Let $C\in {\cal H}$ map a family of sets $\cal B$ into ${\cal B}_0$. If $C$ satisfies either the Tarski axioms (i)$,$ (ii) or (i)$,$ (iii)$,$ or the *-transfer 
*(i)$,$ *(ii)$,$ or *(i)$,$ *(iii) of these axioms$,$ then $C$ is called a {\it subtle consequence operator}.  For example$,$ if $\r C \in {\cal C},$ then  it is immediate that $\Hyper {\b C}\colon \Hyper {(\power {\theta(\r A)}} \to \Hyper {\power {\theta (\r A))}}$ satisfies *(i) and *(ii) for the family of all internal 
subsets of $\Hyper {(\theta (\r A))}.$ This $\Hyper {\b C}$ is a subtle consequence operator. For any set $A \in {\cal N},$ let $\sig A =\{\hyper a\mid a \in A\}.$ (In general, this definition does not correspond to that used by other authors.) If for a subtle consequence operator $C$ there does not exist some similarly defined $D \in {\cal N}$ such that $C = \sig 
{\b D}$ or $C = \Hyper {\b D}$$,$ then $C$ is called a {\it purely}
subtle consequence operator. Let infinite $A \subset {\cal E}$ and $B = \hyper A - \sig A.$ Then the identity $I\colon \power B - \power B$ is a purely subtle consequence operator.\pars
There are certain technical procedures associated with the $\sigma$ map that take on a specific significance for consequence operators. Recall that $\cal N$ is closed under finitely many power set or finite power set 
iterations. Let $X,\ Y \in {\cal N}.$ It is not difficult to show that if ${\cal P}\colon \power X \to Y,$ then for each $A \in \power X,\ \Hyper (\power A) = \Hyper {\cal P}(\hyper A).$ Moreover$,$ if $F\colon \power X \to Y,$ where $F$ is the finite power set operator$,$ then for each $A \in \power A,\ \Hyper {(F(A))}=
\hyper F (\hyper A).$ If $\r C \in \cal C$ and $X \subset \r W$$,$ then $\b C\colon \power X \to \power X$ has the property that for each $A \in \power X,\ \Hyper {(\b C(A))} = \Hyper {\b C}(\hyper A).$\pars
Recall that we identify 
each $\Hyper n \in \hypernat$ with $n \in \nat$ since 
$\Hyper n$ is but a constant sequence with the value $n.$ Utilizing this fact$,$ we have the following
straightforward lemma the proof of which is omitted.\parm
\indent {\bf LEMMA 4.1}.

 (i) {\it Let $A \in \cal N.$ Then 
$^\sigma (F(A)) = F(^\sigma A)$. If also $A \subset ({\cal W} \cup {\cal E}),$ then $^\sigma 
(F(A)) = F(A).$}\pars

(ii) {\it Let ${\rm C} \in {\cal C}^\prime,\ \r B \subset \r X 
\subset \cal W.$} \pars 

(a) {\it $^\sigma (\b C(\b B)) = \b C(\b B).$}\pars 

(b) {\it $\Hyper {\b C}\bigm| \{\hyper {\b A}\mid {\b A} \in \power X \}=
 \{ (\hyper {\b A},\hyper {\b B}) 
\mid (\b A,\b B) \in \b C \} = {^\sigma \b C}.$} \pars 

(c) {\it If $\r F \in \r F(\r B),$ 
then $^\sigma (\b C(\b F)) \subset (^\sigma \b C )(^\sigma {\b F}) = 
(^\sigma \b C)(\b F).$ Also 
$^\sigma (\b C(\b B)) \subset {(^\sigma \b C) (\hyper {\b B})}$ and$,$ in general$,$ 
$^\sigma (\b C(\b B)) \not= (^\sigma \b C)(\hyper {\b B}),\ ^\sigma (\b C(\b F)) \not= 
(^\sigma \b C)(F).$} \pars 

(d) {\it If $C \in {\cal C}_f^\prime,$ then 
$^\sigma (\b C(\b B))= \bigcup \{^\sigma (\b C(\b F)) \mid \b F \in {^\sigma (F(\b B))}\} 
= \bigcup \{ ^\sigma (\b C(\b F)) \mid \b F \in F(^\sigma {\b B})\} = \bigcup \{\b C(\b F) 
\mid \b F \in F(\b B) \}.$}\parm

(A duplicate lemma holds, where $A \in \cal N$ and $C \in \cal C$, where $\cal C$ set of consequence operators defined on subsets of W if W is included as a subset of the ground set. The difference is that the ``bold'' notion does not appear.)\parm
Throughout the remainder of this paper$,$ we remove from $\cal C$ the one and only one inconsistent consequence operator $U.$ Thus notationally we let $\cal C$ denote the set of all consequence operators defined on infinite $ \r L = \r X \subset \r W$ with the exception of $U.$ Two types of chains will be investigated. Let \ T\ be any chain in $\langle {\cal C},\leq \rangle$ and
\  $\r T'$ \ be any chain with the additional property that for each $\r C \in \r T'$ there exists some $\r C' \in \r T'$ such that $\r C < \r C'.$\parm

{\bf THEOREM 4.2} {\it There exists some $C_0 \in \Hyper {\b T}$ such that for each $\r C \in \r T,\ \Hyper {\b C}\leq C_0.$ There exists some $C'_0 \in \Hyper {\b T'}$ such that $C'_0$ is a purely subtle consequence operator and for each $\r C \in \r T',\ \Hyper {\b C} < C_0'.$ Each member of $\Hyper {\b T}$ and $\Hyper {\b T'}$ are subtle consequence operators.}\pars
PROOF. Let $R = \{(x,y)\mid x,\ y \in \b T\}\ {\rm and}\ x\leq y \}$ and 
$R'\{(x,y)\mid x,\ y  \in \b T' \ {\rm and}\ x < y\}.$ In the usual manner$,$ it follows that $R$ and $R'$ are concurrent on the set $\b T$ and $\b T'$ respectively. Thus there is some $C_0 \in \Hyper {\b T}$ and $C'_0\in \Hyper {\b T'}$ such that for each $\r C \in \r T$ and $\r C' \in \r T,\ \Hyper {\b C} \leq C_0$ and $\Hyper {\b C'} < C_0'$ since $\Hyper {\cal M}$ is an enlargement. Note that the members of $\Hyper {\b T}$ and $\Hyper {\b T'}$ are defined on the set of all internal subsets of $\Hyper {\b L}.$
However$,$ if there is some similarly defined $D \in \cal N$ such that 
$C_0$ or $C_0' = \sig D,$ then since $\sig D$ is only defined for *-extensions of the (standard) members of $\power {\r L}$ and each 
$E \in \Hyper {\b T}\ {\rm or}\ \Hyper {\b T'}$ is defined on the internal subsets of $\Hyper {\b L}$ and there are internal subsets of $\Hyper {\b L}$ that are not *-extensions of standard sets we would have a contradiction. Of course$,$ each member of $\Hyper {\b T}$ or $\Hyper {\b T'}$ is a subtle consequence operator. Hence each $E \in \Hyper {\b T}$ or $\Hyper {\b T'}$ is either equal to some $\Hyper {\b C},$ where $\r C \in \r T$ or $\r C\in \r T'$ or it is a purely subtle consequence operator. Now there does not exist a $D \in \cal N$ such that $C_0' = \Hyper {\b D}$ since $C_0' \in \Hyper {\b T'}$ and $\Hyper {\b C} \not= C_0'$ for each $\Hyper {\b C} \in \sig {\b T'}$ would yield the contradiction that $\Hyper D \in \Hyper {\b T'}- \sig {\b T'}$ but $\Hyper {\b D} \in \sig {\bf {\cal C}}.$ Hence $C_0'$ is a purely subtle consequence operator. This completes the proof.\parm
Let $\r C \in \r T'$. Since $\Hyper {\b C} < C_0'$$,$ then $C'_0$ is  ``more
powerful'' than any $\r C\in \r T'$ in the following sense. If $\r B \in \power {\r L},$ then for each $\r C \in \r T'$ it follows that $\b C({\b B}) \subset \Hyper(\b C(\b B)) = \Hyper {\b C}(\Hyper {\b B}) \subset C_0'(\Hyper {\b B}).$ Recall that$,$ for $\r C \in \cal C,$ a set $\r B \in \power {\r L}$ is called a {\it C-deductive system} if $\r C(\r B) = \r B.$ From this point on$,$ all results are restricted to chains in $\langle {\cal C}_f, \leq \rangle.$\pars
{\bf THEOREM 4.3.} {\it Let $\r C \in \r T \cup \r T'$ and $\r B \in \power {\r L}.$ Then there exists a 
*-finite $F_0 \in \Hyper (F(\b B))$ such that $\b C(\b B)\subset \Hyper {\b C}(F_0) \subset \Hyper {\b C}(\Hyper {\b B}) = \Hyper (\b C({\b B}))$ and $\Hyper {\b C}(F_0) \cap \b L = \b C(\b B) = \Hyper {\b C}(F_0) \cap \b C(\b B).$}\pars
PROOF. Consider the binary relation $Q = \{(x,y)\mid x \in \b C(\b B),\ y\in F(\b B)\ {\r and}\ x\in \b C(y)\}.$ By axiom (iii)$,$ the domain of $Q$ is 
$\b C(\b B).$ Let $(x_1,y_1),\ldots,(x_n,y_n) \in Q.$ By theorem 1 in [6$,$ p. 64]$,$ (the monotone theorem) we have that $\b C(y_1)\cup \cdots\cup \b C(y_n) \subset \b C(y_1 \cup \cdots\cup y_n).$ Since $F = y_1 \cup \cdots \cup y_n \in F(\b B),$ then $(x_1,F),\ldots,(x_n,F) \in Q.$ Thus $Q$ is concurrent on ${\b C}(\b B).$ Hence there exists some $F_0 \in \Hyper {(F(\b B))}$ such that $\sig(\b C(\b B)) = \b C(\b B)\subset \Hyper {\b C}(F_0) \subset \Hyper {\b C}(\Hyper {\b B}) = \Hyper (\b C(\b B)).$ Since $\sig {\b L} = \b L,$ then $\Hyper {\b C}(F_0) \cap \b L = 
\b C(\b B) = \Hyper {\b C}(F_0) \cap \b C(\b B).$\parm
{\bf COROLLARY 4.3.1} {\it If $\r C\in {\cal C}_f$ and $\r B \in \power {\r 
L}$ is a C-deductive system$,$ then there exists a *-finite $F_0 \subset \Hyper 
{\b B}$ such  that  $\Hyper {\b C}(F_0) \cap \b L = \b B.$}\pars 
PROOF. Simply consider the one element chain $\r T = \{ \r C\}.$\parm 
{\bf COROLLARY 4.3.2.} {\it Let $\r C \in {\cal C}_f.$ There there exists a *-finite $F_1 \subset \Hyper {\b L}$ such that for each C-deductive system $\r B 
\subset \r L,\ \Hyper {\b C}(F_1)\cap \b B = \b B.$}\pars 
PROOF. Let $\r T = \{\r C\}$ and the  ``B'' in theorem 4.3 equal \ L. The 
result now follows in a straightforward manner.\parm 
{\bf THEOREM 4.4.} {\it Let $\r B \in \power {\r L}.$\pars 
{\rm (i)} There exists a *-finite $F_B \in \Hyper {(F(\b B))}$ and a subtle 
consequence operator $C_B  \in \Hyper {\b T}$ such that for all $\r C \in \r 
T,\ \sig(\b C(\b B)) = \b C(\b B) \subset C_B(F_B).$\pars 
{\rm (ii)} There exists a *-finite $F_B \in \Hyper {(F(\b B))}$ and a purely subtle 
consequence operator $C_B' \in \Hyper {\b T'}$ such that for all $\r C\in \r 
T',\ \sig(\b C(\b B)) = \b C(\b B) \subset C_B'(F_B).$}\pars 
PROOF.  Consider the ``binary'' relation $Q=\{((x,z),(y,w))\mid x \in \b T,\ 
y\in \b T,\ w \in F(\b B),\ z\in x(w),\ z\in \power {\b L}, z \in 
x(w),\ {\rm and}\ x(w) \subset y(w)\}.$ Let $\{((x_1,z_1),(y_1,w_1)),\ldots, 
((x_n,z_n),(y_n,w_n))\} \subset Q.$ Notice that $F = w_1\cup\cdots\cup w_n \in 
F(\b B)$ and for the set
$K =\{x_1,\ldots, x_n\},$ let $D$ be the largest member of $K$ with respect 
to $\leq.$ It follows that $z_i \in x_i(w_i) \subset x_i(F) \subset D(F)$ for 
each $ i = 1,\ldots, n.$ Hence $\{((x_1,z_1),(D,F)),\ldots,((x_n,z_n),(D,F))\} 
\subset Q$ implies that $Q$ is concurrent on its domain. Consequently$,$ there 
exists some $(C_B,F_B) \in \Hyper {\b T} \times \Hyper (F(\b B))$ such that 
for each $(x,z) \in \ {\rm domain\ of}\ Q,\ (\Hyper (x,z), (C_B,F_B))\in 
\Hyper Q.$  Or$,$ for each $(u,v) \in \ \sig ({\rm domain \ of}\ Q),\ ((u,v),(C_B,F_B)) 
\in \Hyper Q.$ Let arbitrary $\r C \in \r T$ and $b \in C(B).$ Then 
there exists some $F' \in F(\b B)$ such that $\b b \in \b C(F').$ Thus $(\Hyper 
{\b C},\Hyper {\b b}) \in \ \sig ({\rm domain\ of}\ Q).$ Consequently$,$ for each 
$\r C \in \r T$ and $\b b \in \b C(\b B),\  \b b =\Hyper {\b b} \in (\Hyper {\b C})(F_B) 
\subset C_B(F_B).$ This all implies that for each $\r C \in \r T,$ 
$\sig (\b C(\b B)) = \b C(\b B) \subset C_B(F_B).$\pars 
(ii) Change the relation $Q$ to $Q'$ by requiring that $x \not= y.$ Replace 
$D$ in the proof of (i) above with $D'$ is greater than and not equal to the 
largest member of $K.$ Such a $D'$ exists in $\b T'$ from the definition of 
$\b T'.$ Continue the proof in the same manner in order to obtain $C_B'$ and 
$F_B'.$ The fact that $C_B'$ is a purely subtle consequence operator follows 
in the same manner as in the proof of theorem 4.2.\parm
{\bf COROLLARY 4.4.1} {\it There exists a \r [resp. purely\r ] subtle consequence 
operator $C_L \in \Hyper {\b T}$ \r [resp. $\Hyper {\b T'}$\r ] and a *-finite $F_L 
\in \Hyper (F(\b L))$ such that for all $\r C \in \r T$ \r [resp. $\r T'$\r ]
and each $\r B \in \power {\r L},\ \b B\subset \b C(\b B) \subset C_L(F_L).$}
\pars
PROOF. Simply let  ``B'' in theorem 4.4 be equal to \ L. Then there exists a 
[resp. purely] subtle $C_L \in \Hyper {\b T}$ [resp. $\Hyper {\b T'}$] and 
$F_L \in \Hyper (F(\b L))$ such that for all $\r C\in \r T$ [resp. $\r T'$] 
$\b C(\b L) \subset C_L(F_L).$ If $\r B \in \power {\r L}$ and $\r C \in \r T$ 
[resp. $\r T'$]$,$ then $\r B \subset \r C(\r B)\subset \r C(\r L).$ Thus for 
each $\r B \in \power {\r L}$ and $\r C\in \r T$ [resp. $\r T'$] 
$\b B \subset \b C(\b B) \subset C_L(F_L)$ and the theorem is 
established.\parm
The nonstandard results in this section have important applications to 
mathematical philosophy. We present two such applications. Let $\cal F$ be the 
symbolic alphabet for any formal language \ L \ with the usual assortment of 
primitive symbols [10$,$ p. 59]. We note that it is possible to mimic the 
construction of \ L \ within $\cal E$ itself. If this is done$,$ then it is not 
necessary to consider the map $\
$ and we may restrict our 
attention entirely to the sets ${\cal E}$ and $\Hyper {\cal E}.$ \pars
Let \ S \ denote the predicate consequence operator by the standard rules for 
predicate (proof-theory) deduction as they appear on pages 59 and 60 of 
reference [10]. Hence $\r A \in \power {\r L},\ \r S(\r A) =\{x\mid x \in \r L\ 
{\rm and}\ \r A \vdash x\}.$ It is not difficult to restrict the modus ponens 
rule of inference in such a manner that a denumerable set $\r T' = \{\r 
C_n\mid n \in \nat \}$ of consequence operators defined on $\power {\r L}$ is 
generated with the following properties.\parm
(i) For each $ \r A \in \power {\r L},\ \r S(\r A)= \bigcup \{ \r C_n(\r 
A)\mid n\in \nat \}$ and $\r C_n \not= \ \r S$ for any $n \in \nat.$\pars
(ii) For each $\r C \in \r T'$ there is a $\r C'$ such that $\r C < \r C'$ [5$,$ 
p.57]. Let $\r A \in \power {\r L}$ be any S-deductive system. The $\r A =\r 
S(\r A)= \bigcup \{\r C_n(\r A)\mid n\in \nat\}$ yields that \ A \ is
a $\r C_n$-deductive system for each $n \in \nat.$ Thus  \ S and \ $\r C_n$  
$n \in \nat$ are consequence operators defined on $\power {\r A}$ as well as
on $\power {\r L}.$\parm 
{\bf THEOREM 4.5.} {\it Let \ $\r L$ \ be a first-order language and $\r A \in 
\power {\r L}.$ Then there exists a purely subtle $C_1 \in \Hyper {\r T'}$ and 
a *-finite $F_1 \in \Hyper (F(\b A))$ such that for each $\r B \in \power {\r 
A}$ and each $\r C \in \r T'$\pars
{\rm (i)} $\b C(\b B) \subset C_1(F_1),$\pars
{\rm (ii)} $\b S(\b B) \subset C_1(F_1) \subset \Hyper {\b S(F_1)} \subset \Hyper 
{(\b S(\b A))}.$\pars
{\rm (iii)} $\Hyper {\b S(F_1)} \cap \b L = \b S(\b A) = C_1(F_1)\cap \b L.$}\pars
PROOF. The same proof as for corollary 4.4.1 yields that there is some purely 
subtle $C_1 \in\Hyper {\b T'}$ and $F_1 \in \Hyper {(F(\b A))}$ such that for 
each $\r B \in \power {\r A}$ and each $\r C \in \r T',\ \b C(\b B) \subset 
C_1(F_1)$ and (i) follows. From (i)$,$ it follows that $\bigcup \{ \sig(\b C
(\b B))\mid \r C\in \r T'\}= \bigcup \{\b C(\b B)\mid \r C  \in \r T'\} 
= \b S(\b B) = \sig(\b S(\b B))\subset C_1(F_1)$ and the first part of (ii) 
holds. By 
*-transfer $C_1 < \Hyper {\b S}$ and $C_1$ and $\Hyper {\b S}$ are defined on 
internal subsets of $\Hyper {\b A}.$ Thus $C_1(F_1) \subset \Hyper {\b S(F_1)}
\subset \Hyper {\b S}(\Hyper {\b A}) = \Hyper (\b S(\b A))$ by the *-monotone 
property. This completes (ii). Since $\b S(\b A) \subset C_1(F_1) \subset 
\Hyper {\b S}(F_1) \subset \Hyper (\b S(\b A))$ from (ii)$,$ then (iii) follows
and the theorem is proved.\parm
REMARK 4.6. Of course$,$ it is well known that there exists some $F \in \Hyper 
(F(\b A))$ such that $\b S(\b A) \supset \b A \subset F \subset \Hyper {\b 
A}$ and *-transfer of axiom (i) yields that $\Hyper {\b S(F)} \subset \Hyper 
{\b S}(\Hyper {\b A}) = \Hyper {(\b S(\b A))}.$ However$,$ $F_1$ of theorem 4.5 
is of a special nature in that the purely subtle $C_1$ applied to $F_1$ yields 
the indicated properties. Also theorem 4.5 holds for many other infinite 
languages and deductive processes.\pars
Let \ L \ be a language and let $M$ be a structure in which \ L \ can be 
interpreted in the usual manner. A consequence operator \ C \ is  {\it sound} 
for $M$ if whenever $\r A \in \power {\r L}$ has the property that 
$M\models \r A,$ then $M\models \r C(\r A).$ As usual$,$ $T(M) =\{ x\mid
x \in \r L \ {\rm and}\ M \models x\}.$ Obviously$,$ if \ C \ is sound for $M$$,$ 
then $T(M)$ is a C-deductive system. \pars
Corollary 4.3.1 implies that there exists *-finite $F_0 \subset \Hyper 
{(\b T(M))}$ such that $\Hyper {\b C}(F_0)\cap \b L = \b T(M).$ Notice that the fact 
that $F_0$ is *-finite implies that $F_0$ is *-recursive. Moreover$,$ trivially$,$ $F_0$ is a 
*-axiom system for $\Hyper {\b C(F_0)},$ and we do not lack knowledge about 
the behavior of $F_0$ since any formal property about $\b C$ or recursive sets$,$ 
among others$,$ must hold for $\Hyper {\b C}$ or $F_0$ when property 
interpreted.
If \ L \ is a first-order language$,$ then \ S \  is sound for first-order 
structures. Theorem 4.5 not only yields a *-finite $F_1$ but a purely subtle 
consequence operator $C_1$ such that$,$ trivially$,$ $F_1$ is a *-axiom for 
$C_1(F_1)$ and for $\Hyper {\b S(F_1)}.$ In this case$,$ we have that 
$\Hyper {\b S}(F_1)\cap \b L = \b T(M) = C_1(F_1)\cap \b L.$ By the use of 
internal and external objects$,$ the nonstandard logics $\{C_1,\Hyper {\b L}\}$ and $\{\Hyper {\b S}, \Hyper {\b L}\}$ 
technically by-pass a portion of G\"odel's first incompleteness theorem.\pars

By definition $\r b \in \r S(\r B),\ \r B \in \power {\r L}$ iff there is a 
finite length ``proof'' of \ b \ from the premises \ B. It follows$,$ that for 
each $b \in \Hyper {\b S}(F_1)$ there exists a *-finite length proof of $b$ from 
a *-finite set of premises $F_1.$  In particular, for each $b \in \b T(M).$ From the first line of this remark, where $\b A = \b T(M),$ there is a  *-finite $F_2 \subset \Hyper {\b L}$ such that $\sig{\b T(M)}=  \b T(M) \subset F_2$, and, in this case, $\b T(M) \subset \Hyper {\b S}(F_2).$ \bigskip
\centerline{\bf References}\parm
\noindent 1. \ A. Abian: The theory of sets and transfinite arithmetic, W. B. Saunders Co., Philadelphia and London, 1965.\pars
\noindent 2. \ W. Dziobiak: The lattice of strengthenings of a strongly finite 
consequence operator$,$ Studia Logica$,$ {\bf 40}(2) (1981)$,$ 177-193.\pars
\noindent 3. \ J. R. Geiser: Nonstandard logics$,$ J. Symbolic Logic$,$ {\bf 
33}(1968)$,$ 236-250.\pars
\noindent 4. \ C. W. Henson: The isomorphism property in nonstandard analysis 
and its use in the theory of Banach spaces$,$ J. Symbolic Logic$,$ {\bf 39}(1974)$,$ 
717-731.\pars
\noindent 5.\ R. A. Herrmann: The mathematics for mathematical philosophy$,$ 
Monograph \#130$,$ Institute for Mathematical Philosophy Press$,$ Annapolis$,$ 
1983. (Incorporated into ``The Theory of Ultralogics$,$''  http://xxx.arxiv.org/abs/math.GM/9903081 and 9903082.)\pars
\noindent 6. \ R. A. Herrmann: Mathematical Philosophy$,$ Abstracts A. M. S. 
{\bf 2}(6)(1981)$,$ 527.\pars
\noindent 7. \ T. Tech: The axiom of choice, North-Holland and American Elsevier Publishing Co., 1973.\pars
\noindent 8. \ J. Los and R. Suszko: Remarks on sentential logics$,$ Indag. 
Math.$,$ {\bf 20}(1958)$,$ 177-183.\pars
\noindent 9. \ A. A. Markov: Theory of algorithms$,$ Amer. Math. Soc. 
Translations$,$ Ser. 2$,$ {\bf 15}(1960)$,$ 1-14.\pars
\noindent 10. \ E. Mendelson: Introduction to mathematical logic$,$ D. Van 
Nostrand$,$ New York$,$ 1979.\pars
\noindent 11. \ A. Robinson: On languages which are based on non-standard 
arithmetic$,$ Nagoya Math. J.$,$ {\bf 22}(1963)$,$ 83-118.\pars
\noindent 12. A. Tarski: Logic$,$ semantics and metamathematics$,$ Oxford 
University Press$,$ Oxford 1956.\pars
\noindent 13. \  R. W\'ojcicki: Some remarks on the consequence operation in 
sentential logics$,$ Fund. Math. {\bf 68}(1970)$,$ 269-279.\pars
\bigskip\bigskip
\noindent{Robert A. Herrmann}\par
\noindent{Mathematics Department$,$}\par
\noindent{U. S. Naval Academy$,$}\par
\noindent{ANNAPOLIS$,$ MD 21402}\par
\noindent{U. S. A.}\par
\end